\documentclass[11pt]{article}
\usepackage{amssymb}

\newtheorem{theorem}{Theorem}[section]
\newtheorem{prop}{Proposition}[section]
\newtheorem{lemma}[theorem]{Lemma}

\newtheorem{conj}[theorem]{Conjecture}

\def\Gal{\mathop{\rm Gal}\nolimits}
\def\x{$\hfill\rlap{$\sqcup$}\sqcap$\bigskip}

\long\def\atn#1{\setbox10=\vbox{#1}\par
\hbox{\vbox{\hbox to 0pt{\hss\vrule width 1pt height \ht10
\ \vrule width 1pt height \ht10\kern 5pt}}\box10}}

\def\eb#1{{\color{blue}#1}}
\def\eb#1{{#1}}

\let\lfc=\longmapsto

\let\cj=\overline

\let\ms=\medskip

\def\Z{{\msb Z}}

\def\F*2g{{\msb F}^*_{2^g}}
\def\f#1{{\msb F}_{#1}}

\def\F#1{{\msb F}_{2^{#1}}}

\def\pt#1{\left(#1\right)}

\def\biq#1{\quad \hbox{ #1 }\quad }

\def\no{n$^\circ$}
\newcount\annee\annee=\year\advance\annee by -2000
\def\og{\leavevmode\raise.3ex\hbox{$\scriptscriptstyle\langle\!\langle$}}
\def\fg{\leavevmode\raise.3ex\hbox{$\scriptscriptstyle\rangle\!\rangle$}}

\font\tenmsb=msbm10
\font\sevenmsb=msbm7
\font\fivemsb=msbm5
\newfam\msbfam%Family 9
\def\msb{\fam\msbfam\tenmsb}%
\textfont\msbfam=\tenmsb \scriptfont\msbfam=\sevenmsb%
\scriptscriptfont\msbfam=\fivemsb%

\begin{document}

\title{Functions of degree $4e$\\ that are not APN Infinitely Often}
\author{Fran\c cois Rodier}
\date{}
\maketitle

\begin{abstract}
We prove a necessary condition for some polynomials of degree $4e$ ($e$ an odd number) to be   APN over $\mathbb{F}_{q^n}$ for large $n$, and we investigate  the polynomials $f$ of degree 12.
\end{abstract}


\section{Introduction}
\ms

The vector Boolean functions are used in cryptography to
construct block ciphers and an important criterion on these
functions  is their  high resistance to differential cryptanalysis.

Let  $q=2^n$ for some positive integer $n$.
A function $f~:~\f q \longrightarrow \f q$ is said to be \emph{almost perfect nonlinear} (APN) on $\f q$
if the number of solutions in $\f q$ of the equation
$$f(x+a)+f(x)=b$$
is at most 2, for all $a,b\in \f q$, $a \not=0$.
Because $\f q$ has characteristic 2, the number of solutions
to the above equation must be an even number, for any function $f$ on $\f q$.
This kind of function has a good resistance to differential cryptanalysis as was  proved by  Nyberg in \cite{ny}.

So far, the study of APN functions has focused on power functions.
Recently it was generalized to other functions, particularly
polynomials (Carlet, Pott and al. \cite{bcfl, ekp, ep}) or polynomials
on small fields (Dillon \cite{JD}).
On the other hand, several authors (Berger, Canteaut, Charpin, Laigle-Chapuy \cite{bccl}, Byrne and
McGuire \cite{bg} or Jedlicka \cite{je}) showed that APN functions did not exist in certain cases.

\eb{We would like to have a complete classification of all APN functions. Indeed there are many classes of function for which it can be shown
that each function is APN for at most a finite number of extensions.
 So we fix a finite field $\f{q}$ and a function $f:\f{q}\to \f{q}$ given by a polynomial in $\f{q}[x]$ and
we set the question of whether this function can be APN  for an infinite number of extensions of $\f{q}$.}

In this approach, Hernando and McGuire  \cite{HM} showed a result on the classification of APN monomials which has been conjectured for 40 years: the only exponents such that the monomial $x^d$ are APN over  infinitely many
extension of $\f2$ are of the form $2^i+1$ or $4^i-2^i+1$. One calls these exponents {\sl exceptional exponents}.
Then it is natural to
 formulate for polynomial functions
the following conjecture.
\begin{conj}[Aubry, McGuire and Rodier]
A polynomial on $\f q$\break can
be APN for an infinity of extensions of $ \f q $ only if it is CCZ equivalent (as was defined by Carlet, Charpin and   Zinoviev in \cite{ccz}) to a monomial $ x ^ t $
where $ t $ is an exceptional exponent.
\end{conj}
A means to prove this conjecture is to remark that the APN property is equivalent to the fact that the rational points of a certain algebraic surface $X$ in a 3-dimensional space linked with the polynomial $f$ defining the  Boolean function are all in a surface $V$ made of 3 planes and independent of $f$.
We define the surface $X$ in the three dimensional affine space $ \mathbb{A}^3$ by
\[
\phi (x,y,z) = \frac{f(x)+f(y)+f(z)+f(x+y+z)}{(x+y)(x+z)(y+z)}
\]
which is a polynomial in $\mathbb{F}_q [x,y,z]$.
When this surface is irreducible (or when it has an irreducible component defined over the field of definition of $f$),
a Weil's type bound may be used to approximate the number of rational points of this surface. When it is too large the surface is too big to be contained in the surface $V$, and the function $f$ cannot be APN.

A relation between the degree and the number of variables on a large  number of  Boolean function was obtained by this means by the author \cite{FR}.  This enabled him, with Aubry and McGuire to prove the conjecture  in a number  of cases.

We begin by recalling some known results. Then we give some new results for polynomials of degree $4e$ where $e$ is odd and we investigate completely the case of polynomials of degree 12.

\section{The State of the Art}

The best known examples of APN functions are the Gold functions $x^{2^k+1}$ and
the Kasami-Welch functions $x^{4^k-2^k+1}$.
These functions are defined over $\mathbb{F}_{2}$, and are APN on
any field $\mathbb{F}_{2^m}$ where $gcd(k,m)=1$.

\begin{theorem}[Aubry, McGuire and Rodier, \cite{amr}]
If the degree of the polynomial function $f$ is odd and not a Gold or a Kasami-Welch number then $f$ is not APN over $\mathbb{F}_{q^n}$ for all $n$ sufficiently large.
\end{theorem}

In the even degree case, we can state the result when half of the degree is odd, with an extra minor condition.

\begin{theorem}[Aubry, McGuire and Rodier, \cite{amr}]
\label{degre2e}
If the degree of the polynomial  function $f$ is $2e$
with $e$ odd,   and if $f$ contains a term of odd degree, then $f$ is not APN over $\mathbb{F}_{q^n}$ for all $n$ sufficiently large.
\end{theorem}

We have some results on the polynomials of Gold degree $d=2^k+1$.

\begin{theorem}[Aubry, McGuire and Rodier, \cite{amr}]
Suppose $f(x)=x^d+g(x)$ where $\deg (g) \leq 2^{k-1}+1$ .
Let $g(x)=\sum_{j=0}^{2^{k-1}+1} a_j x^j$.
Suppose moreover that there exists a nonzero coefficient $a_j$ of $g$ such that
$\phi_{j} (x,y,z)$ is absolutely irreducible (where $\phi_i (x,y,z)$ denote the polynomial $\phi (x,y,z)$ associated to $x^i$).
Then   $f$ is not APN over $\mathbb{F}_{q^n}$ for all $n$ sufficiently large.
\end{theorem}

\section{New results: polynomials of degree $4e$ where $e$ is odd}

We have been interested in the polynomials of degrees of the form $4e$, where $e$ is an odd number.
This case is far more intricate than the previous cases, because there are some polynomials which are CCZ equivalent to monomials.

Here are some partial results.

\begin{theorem}
\label{geom}
If the degree of the polynomial function $f$ is even such that
$\deg(f)=4e$ with $e\equiv3\pmod4$,
and if the polynomials of the form
$$(x+y)(y+z)(z+x)+P$$ with
\begin{equation}
\label{poly}
P(x,y,z)= c_1 (x^2 + y^2 + z^2) + c_4 (x y+   x z+  z y)+ b_1 (x+  y+   z)+d
\end{equation}
for $c_1, c_4, b_1, d\in \f{q^3}$,
do not divide $\phi$
then $f$ is not APN over $\mathbb{F}_{q^n}$ for $n$ large.
\end{theorem}

We have more precise results for polynomials of degree 12.
\begin{theorem}
\label{d12}
If the degree of the polynomial $f$ defined over $\f q$ is 12,    then either $f$ is not APN over $\mathbb{F}_{q^n}$ for large $n$ or $f$ is CCZ equivalent to the Gold function $x^3$.
\eb{In this case $f$ is of the form
$$L(x^3)+L_1 \biq{or} (L(x))^3+L_1$$
where $L$ is a linearized polynomial
$$x^4+x^2(c^{1+q}+c^{1+q^2}+c^{q+q^2})+x c^{1+q+q^2},$$
$c$ is an element of $\f{q^3}$ such that $c+c^q+c^{q^2}=0$ and $L_1$ is a q-affine polynomial of degree at most 8 (that is a polynomial   whose monomials are of degree  0 or a power of 2).}
\end{theorem}

\section{Preliminaries}

We first eliminate some trivial functions.
The following propositions are easy to prove.

\begin{prop}
\label{p2}
The class of APN functions is invariant by adding a   $ q $-affine polynomial.
\end{prop}

\begin{prop}
The kernel of the map
\[
f \lfc \frac{f(x)+f(y)+f(z)+f(x+y+z)}{(x+y)(x+z)(y+z)}
\]
is made of   $ q $-affine polynomials.
\end{prop}

\begin{prop}
\label{irrcomp}
When the surface  has an absolutely irreducible component defined over the field of definition of $f$ which is not one of the planes $(x+y)(y+z)(z+x)=0$,  the function $f$ cannot be APN for infinitely many extension of $\f q$.
\end{prop}

This proposition is a consequence of Theorem 4.1 and Corollary 3.1 of \cite{FR}.

\section{Proof of Theorem \ref{geom}}

We suppose that  the degree of the polynomial function $f$ is
$4e$ with $e\equiv3\pmod4$, and that $f$   is APN for infinitely many extensions.
We want to show that $\phi$ is divisible by polynomials of the form
$(x+y)(x+z)(y+z)+P$ with $P$ as in (\ref{poly}).

\subsection{Some lemmas}

Before proving Theorem \ref{geom} that let us state a few lemmas.
Let  $\cj X$ be the projective closure of $X$.

\begin{lemma}\label{irred}
Let $H$ be a projective hypersurface.
If $\overline X \cap H$ has a reduced absolutely irreducible component defined over $\mathbb{F}_q$
then
$\overline X$ has an absolutely irreducible component defined over $\mathbb{F}_q$.
\end{lemma}

This is Lemma 2.1 in \cite{amr}.

\begin{lemma}
\label{phi_d}
Suppose $d$ is even and write $d=2^je$ where $e$ is odd.
In  $\overline X\cap H$ we have
\[
\phi_d (x,y,z)=\phi_e(x,y,z)^{2^j} ((x+y)(x+z)(y+z))^{2^j-1}
\]
\end{lemma}

See the proof of Lemma 2.2 in \cite{amr}.

\subsection{The component $X_0$}
\label{X_0}

As a consequence of Proposition \ref{irrcomp}
 no absolutely irreducible component of $X$ is defined over $\f  q$, except perhaps $x+y=0$, $z+y=0$ or $x+z=0$.

If some component of $X$ is equal to one of these planes, then by symmetry of $\phi$ in $x$, $y$ and $z$, all of them are component of $X$, which implies that $(x+y)(z+y)(x+z)$ divides $\phi$. Let us suppose from now on that this is not the case.

%Sketch of proof of Theorem \ref{geom}.
Let   $H_{\infty}$ be the plane at infinity of the space  $ \mathbb{A}^3$  and  $X_{\infty}=H_{\infty}\cap \cj X$.
The equation of $X_{\infty}$ is $\phi_d (x,y,z)=0$ where $\phi_d$ is the $\phi$ function associated to the monomial $x^d$.
As we have in this case (cf.  Lemma  \ref{phi_d}):
$$\phi_d (x,y,z)=\phi_e(x,y,z)^{4} ((x+y)(x+z)(y+z))^{3}$$
the curve $X_{\infty}$ is not reduced. Note that since $e\equiv3\pmod4$, the curve $\phi_e(x,y,z)=0$ is irreducible (cf. \cite{JW}).

Let $X_0$ a reduced absolutely irreducible component of $\cj X$ which contains the line $x+y=0$ in $H_{\infty}$.

\subsection{Case where $X_0$ contains three copies of the line $x+y=0$ in $H_{\infty}$.}

If $X_0$ contains three copies of this line, then $X_0$ is defined over $\mathbb{F}_q$ by a Galois argument as in the proof of Lemma \ref{irred}.
Indeed, suppose for the sake of contradiction that $X_0$ is not defined over  $\mathbb{F}_q$.
Then $X_0$ is defined over  $\mathbb{F}_{q^t}$ for some $t$.
Let $\sigma$ be a generator for the Galois group $\Gal(\mathbb{F}_{q^t}/\mathbb{F}_{q})$ of $\mathbb{F}_{q^t}$ over $\mathbb{F}_{q}$.
Then $\sigma(X_0)$ is an absolutely irreducible component of $\overline X$ that is
distinct from $X_0$ hence disjoint from $X_0$.
So
$X_0\cap H_\infty$ which contains  three copies of the line  $x+y=0$ in $H_{\infty}$ is disjoint from $\sigma(X_0\cap H_\infty)$
 which contains  three copies of the same line because this line is defined on $\mathbb{F}_2$. This is a contradiction with the fact that
$\overline X\cap H_\infty$ contains only 3 copies of this line.

\subsection{Case where $X_0$ contains two copies of the line $x+y=0$ in $H_{\infty}$.}

If it contains this line exactly two times, then another component will contain one copy of this line and hence will be defined over $\mathbb{F}_q$ by the same argument as before.

\subsection{Case where $X_0$ contains one copy of the line $x+y=0$ in $H_{\infty}$.}

Suppose now that it contains just one time this line. We can suppose that there are two other components  $X_1$ and $X_2$ which contain this line. The component $X_0$ is defined over an extension, say $\mathbb{F}_{q^t}$ of $\mathbb{F}_q$ and we choose for $t$ the smallest possible. Let $G$ be the Galois group of $\mathbb{F}_{q^t}/\mathbb{F}_q$; since $G$ fixes the line $x+y=0$ in $H_\infty$, the group $G$ acts on the $X_i$'s and let us consider the orbit of $X_0$ under this action. If it contains just $X_0$, then $X_0$ is defined over $\mathbb{F}_q$. If it contains $X_0$ and $X_1$ then $G$ fixes $X_2$ and $X_2$ is then defined over $\mathbb{F}_q$. Finally we suppose that it contains the three components.
Then $G$ acts transitively on these 3 components. Let $G_1$ the stabilizer of $X_0$. Then the group $G/G_1$ is isomorphic to $\Z/3\Z$,
and $G_1$ is the only subgroup of $G$ of index 3.

The same is true for the lines $y+z$ and $z+x$.  Hence one gets the same subgroup $ G_1 $.

\subsection{Case where $X_0$ is of degree 1}
\label{plan}

 If $X_0$ is of degree 1, as the intersection with $H_\infty$ contains the line $x+y=0$, the equation of $X_0$ would be $x+y+b=0$ with $b\in\f {q^t}$ and $b\notin\f q$.  In this case   $x+y+b$ would divide $f(x)+f(y)+f(z)+f(x +y+z)$.

As $b\notin \mathbb{F}_q$, by the  action of $G$, $x+y+\rho(b)$ where $\rho\in G$ would be a distinct plane containing the line $x+y=0$ in $H_\infty$. As there are only three distinct component of $X$ containing  the line $x+y=0$ in $H_\infty$ and as $t$ is minimal, this implies that $t=3$.

\eb{
By symmetry of the variables $x, y ,z$ in the expression of $f(x)+f(y)+f(z)+f(x +y+z)$,  $z+y+b$ and $x+z+b$ divide also  $f(x)+f(y)+f(z)+f(x +y+z)$.}
Finally $f(x)+f(y)+f(z)+f(x +y+z)$ is divisible by
$$ (x+y+b)( z+y+b)( x+z+b)
= (x+y)(y+z)(z+x)+b(x^2 + y^2 + z^2 x y+   x z+  z y)+b^3,$$
which is in the form specified by the Theorem.

\subsection{Case where $X_0$  is of degree at least  2.}
If $X_0$ is not a plane it is of degree at least  2.
If $X_0$ does not contain the  curve $\phi_e(x,y,z)=0$ in $H_\infty$, then $X_0$ contains
one of the two lines $(x+z)=0$ or $(y+z)=0$ and $X_1$ and $X_2$ contain the same, by the action of $G$.

Suppose $X_0$   contains $(y+z)=0$. Since $\phi$ is symmetric in $x$, $y$ and $z$, there exists also 3 components $Y_0$, $Y_1$ and $Y_2$ which contains the two lines  $(y+z )=0$ and $(z+x)=0$. Since $\overline X\cap H_\infty$ contains only 3 copies of $(y+z )=0$, this implies that $Y_i=X_i$ up to the order of indices and thus $X_0$ contains also the line $z+x=0$.

Finally, either $X_0$
contains $\phi_e(x,y,z)=0$, or $X_0$,  $X_1$ and $X_2$ contain the three lines $(x+y)(x+z)(y+z)=0$. In this case, the components $X_0$, $X_1$ and $X_2$ are of degree 3.

\subsection{Case where $X_0$ contains the curve $\phi_e(x,y,z)=0$}

If $X_0$ contains the curve $\phi_e(x,y,z)=0$,  it is the same for $X_1$ and $X_2$, therefore $X_0$ can contain only one copy of this curve.
Another component $X_3$ will therefore contain the 4$^{\hbox{\tiny  th}}$, so it will be defined on $\f q$.

\subsection{Case where  $X_0$ is of degree 3 and contains the lines\break  $(x+y)(x+z)(y+z)=0$  in the plane at infinity.}
\label{X_0irred}

So the only problem is when $X_0$ is of degree 3 and contains the lines\break $(x+y)(x+z)(y+z)=0$ at infinity.
The equation of such a surface is of the form $(x+y)(x+z)(y+z)+P(x,y,z)$ where $P$ is a polynomial  of degree at most 2.

Let $\rho$ a generator of $G$.
The equation of $X_1$ is (say) $(x+y)(x+z)(y+z)+\rho(P)(x,y,z)$ and the equation of $X_2$ is   $(x+y)(x+z)(y+z)+\rho^2(P)(x,y,z)$.
Since these polynomials are irreducible (since we have supposed that $X_0$ is irreducible) and  {distinct}, they are prime with each other, therefore
$f(x_0)+f(x_1)+f(x_2)+f(x_0+x_1+x_2)$ is divisible by
\begin{equation}
\label{prodX}
\prod_{i=0}^3\pt{(x+y)(x+z)(y+z)+\rho^i(P)(x,y,z)}
\end{equation}

The equation of the curve $X_\infty$ is
$$\phi_e^4(x,y,z)\prod_{i=0}^3\pt{(x+y)(x+z)(y+z)}$$
so we find that the product (\ref{prodX}) can contain only three summands, hence $\rho^3(P)=P$.
Hence $P$ is defined on $\f{q^3}$ and $X_0$ also.

The product (\ref{prodX}) must be symmetric  in the  variables $x$, $y$, $z$, since if it were not, the image of the product (\ref{prodX}) by some element of the  symmetry group $\cal G$ of the 3 variables would be different, and also divide $f(x_0)+f(x_1)+f(x_2)+f(x_0+x_1+x_2)$, therefore forcing the curve $X_\infty$ to contain more than 3 times the line $x+y=0$.

If $P$ is not symmetric in the  variables $x$, $y$, $z$, then  the orbit of $P$ by the symmetry group $\cal G$ of the 3 variables would be contained in the set
$\{P, \rho(P),\rho^2(P)\}$ since the product (\ref{prodX}) is symmetric.
The orbit of $P$ under $\cal G$ is not reduced to $\{P\}$ since $P$ is not symmetric.
It is not either reduced to 2 elements, because the third element would be symmetric, so it is equal to the set $\{P, \rho(P),\rho^2(P)\}$. The stabilizer of $P$ in $\cal G$ would then be reduced to a transposition.
But the stabilizer of $\rho(P)$ would contain a conjugate transposition, and this transposition would also fix  $P$, as the action of $G$ and $\cal G$ commute. So it is impossible, which proves that $P$ must be symmetric.

Therefore $P$ is of the form
$$P(x,y,z)= c_1 (x^2 + y^2 + z^2) + c_4 (x y+   x z+  z y)+ b_1 (x+  y+   z)+d.$$

\section{Proof of Theorem \ref{d12}}

This theorem is a consequence of  the following proposition.

We put from now on
$A=(x+y)(y+z)(z+x)$.

\begin{prop}
\label{poly12}
The polynomials $\phi$ for $f$ defined over $\f q$ such that $\deg(f)=12$ and that $f$ is APN for an infinite number of extension of $\f q$ are  of the form
\begin{itemize}
\item
$$A^3+\beta A\mu^2+ (A^2+\mu^3)\gamma +A(\gamma^2+\beta ^3)+\beta ^2\gamma\mu+\gamma^3$$
with
$c$ in $\f{q^3}$ such that $c+c^q+c^{q^2}=0$,
$\beta =c^{q+q^2}+c^{1+q^2}+c^{1+q}$,   $\mu=x^2 + y^2 + z^2+x y+   x z+  z y$
and
$\gamma=c^{1+q+q^2}$
\item{\eb{or
$$A^3+\beta A+\gamma^3.$$
}}
\end{itemize}
\end{prop}

We decompose the polynomial  $\phi$ into homogeneous components:
$$\sum_{i=3}^{12} a_i \phi_i(x,y,z)=\phi(x,y,z)= \frac{f(x)+f(y)+f(z)+f(x+y+z)}{(x+y)(x+z)(y+z)}.
$$

\subsection{One more lemma}

Before  proving Proposition \ref{poly12} let us state one more lemma.

\begin{lemma}
\label{divisephi}
The function $x+y$ does not divide $\phi_r(x,y,z)$ for $r$ an odd integer.
\end{lemma}

The function $x+y$ doesn't divide $\phi(x,y,z)$ if and only if the function $(x+y)^2$ doesn't divide $f(x)+f(y)+f(z)+f(x+y+z)$.

 We show easily that   $(x+y)^2$ doesn't divide $x^r+y^r+z^r+(x+y+z)^r$ by using the change of variables $s=x+y$ which gives:
$$x^r+y^r+z^r+(x+y+z)^r=s(x^{r-1}+z^{r-1})+s^2P$$
where $P$ is a polynomial.
\x

\subsection{Proof of Proposition \ref{poly12}}

If the polynomial $f$ is APN for  an infinite number of extension of $\f q$, the polynomial $\phi$ is divisible by $A+P$ for some $P$ by theorem \ref{geom}.

\subsubsection{Reducibility of $A+P$}

In the case where  $A+P$ is reducible one has to consider two cases: either $P=0$ or $P\ne0$.

In the case $P=0$, the polynomial  $\phi$ is divisible by $A$, and so are the homogeneous factors $\phi_r$, so the coefficients $a_r$ are zero if $r$ is odd by the previous Lemma \ref{divisephi}. The polynomial $f$ is therefore the square of a polynomial $f_1$ of degree 6. By Proposition 5.3 of \cite{FR} or Theorem \ref{degre2e} and Proposition \ref{p2} in this paper, it may be APN for an infinity of extension of $\f q$ if and only if it is CCZ equivalent to the function $x^3$.

In the other case, the polynomial $A+P$ has a degree 1 factor which is of the shape $x+y+b$, up to permutation of variables. Then, by section \ref{X_0}, since $P\ne0$, one has $b\ne0$ and even $b\notin\f q$ by Proposition \ref{irrcomp}. As in section \ref{plan}
the polynomials $x+y+\rho(b)$, $x+z+\rho(b)$ and $z+y+\rho(b)$ for $\rho$ in $G$ divide $\phi$ and are prime with each other for $b\notin\f q$, so their product divides $\phi$ and is equal to  $\prod_{i=1}^3 (A+\rho^i(P))$.

If  $A+P$ is irreducible,  the  polynomial
$\prod_{i=1}^3 (A+\rho^i(P))$
divides $\phi$ since it is a product of 3 distinct irreducible polynomials which divide $\phi$.

\subsubsection{The terms of degree 8.}

The polynomial $A+P$ is equal to
$$(x+y)(y+z)(z+x)+  c_1 (x^2 + y^2 + z^2) + c_4 (x y+   x z+  z y)+ b_1 (x+  y+   z)+d.$$
The following polynomial
\begin{equation}
\label{poltotal}
\prod_{\rho\in G} ((x+y)(y+z)(z+x)+\rho( c_1) (x^2 + y^2 + z^2)+\rho( c_4)(x y+   x z+  z y)
+ \rho(b_1)( x+  y+   z) +\rho(d))
\end{equation}
divides $\phi$, so it is equal to $\phi$ since it has the same degree.

The terms of degree 8 in this polynomial are
$$\Big((x+y)(y+z)(z+x)\Big)^2 \Big((x^2+y^2+z^2) (c_1 + \rho(c_1) + \rho^2(c_1)) + (x y+   x z+  z y)(c_4 + \rho(c_4) + \rho^2(c_4))\Big).$$
They must be equal to
$a_{11}\phi_{11}$.
From the lemma \ref{divisephi}, $((x+y)(y+z)(z+x))^2$ does not divides $\phi_{11}$. Hence $a_{11}=0$.
But this implies
$$0=(x^2+y^2+z^2)(c_1 + \rho(c_1) + \rho^2(c_1)) +(x y+   x z+  z y) (c_4 + \rho(c_4) + \rho^2(c_4).$$
Thus
$$ c_1 + \rho(c_1) + \rho^2(c_1)=c_4 + \rho(c_4) + \rho^2(c_4)=0.$$

 \subsubsection{The terms of degree 7.}

The terms of degree 7 in the product
 (\ref{poltotal})
are the term
 $$(b_1+\rho(b_1)+\rho^2(b_1))((x+y)(y+z)(z+x))^2(x+  y+   z)$$ and the terms
 $$\displaylines{
(x+y)(y+z)(z+x)(a_0(x^2 + y^2 + z^2)^2+b_0(x y+   x z+  z y)^2+
 \hfill\cr\hfill
 c_0(x^2 + y^2 + z^2)(x y+   x z+  z y))
}$$
 where $a_0$, $b_0$ and $c_0$ are scalars.
 The term of degree 7 should be equal up to a constant to  $\phi_{10}$ that is $(x+y)(y+z)(z+x)\phi_5^2$.
 Dividing  all terms by $ (x + y) (y + z) (z + x) $, we find that the term
  $$\displaylines{
(b_1+\rho(b_1)+\rho^2(b_1))(x+y)(y+z)(z+x) (x+  y+   z)+a_0(x^2 + y^2 + z^2)^2+
 \hfill\cr\hfill
b_0(x y+   x z+  z y)^2+c_0(x^2 + y^2 + z^2)(x y+   x z+  z y)
}$$
 should be  equal up to a constant to   $\phi_5^2$.
 Now  the term
 $$(x^2 + y^2 + z^2)(x y+   x z+  z y)=x^3 y + x y^3 + x^3 z + x^2 y z + x y^2 z + y^3 z + x y z^2 + x z^3 +  y z^3$$
 contains monomials ($ x^2 y z $) which do not appear in
 $$(x+y)(y+z)(z+x) (x+  y+   z)=x^3 y + x y^3 + x^3 z + y^3 z + x z^3 + y z^3$$
 nor in the square terms
 $(x^2 + y^2 + z^2)^2$ or $(x y+   x z+  z y)^2$. So $c_0=0$.
 Now
 $(b_1+\rho(b_1)+\rho^2(b_1))(x+y)(y+z)(z+x) (x+  y+   z)$  contains monomials $(b_1+\rho(b_1)+\rho^2(b_1))x^3y$  which do not appear in  the square terms
 $(x^2 + y^2 + z^2)^2$ or $(x y+   x z+  z y)^2$  the same way,  whence
 $$b_1+\rho(b_1)+\rho^2(b_1)=0.$$

\subsubsection{The terms of degree 6.}

The polynomial $A+P$ is equal to
\begin{eqnarray*}
&&(x+y)(y+z)(z+x)+  c_1 (x^2 + y^2 + z^2) + c_4 (x y+   x z+  z y)+ b_1 (x+  y+   z)+d \\
&&=(x+y)(y+z)(z+x)+  P_1(x,y,z)+ b_1 (x+  y+   z)+d
\end{eqnarray*}
where $P_1= c_1 (x^2 + y^2 + z^2) + c_4 (x y+   x z+  z y)$.

The following polynomial
$$\prod_{\rho\in G} ((x+y)(y+z)(z+x)+\rho( P_1)
+ \rho(b_1)( x+  y+   z))+\rho(d)$$
 is equal to $\phi$.

The terms of degree 6 in the above expression are:
$$P_1\rho( P_1)\rho^2( P_1)(x,y,z)+P_0(x,y,z)((x+y)(y+z)(z+x))$$
where $P_0$ is a polynomial of degree 3.
This expression must be equal to  $a_9 \phi_9$.
Therefore the polynomial $ x +  y$ must divide the polynomial
$$a_9 \phi_9+P_1\rho( P_1)\rho^2( P_1)(x,y,z).$$
One checks that
$$\phi_9\equiv (x+z)^6 \pmod{x+y}$$
and
$$P_1\equiv c_1 z^2+c_4x^2 \pmod{x+y}.$$
Therefore one has
$$a_9  (x+z)^6= (c_1 z^2+c_4x^2)(\rho(c_1) z^2+\rho(c_4) x^2)(\rho^2(c_1) z^2+\rho^2(c_4) x^2).$$
So
$ x +  z$ divides the second member, thus divides one of the factors, which implies
$c_1=c_4$.

 \subsubsection{Other conditions from the terms of degree 5 and 6.}

If   $c_1 = c_4$ the polynomial $A+P$ becomes
$$(x+y)(y+z)(z+x)+ c_1 (x^2 + y^2 + z^2+x y+   x z+  z y)+b_1(x+y+z)+d.$$
Hence the polynomial (\ref{poltotal}) becomes
\begin{equation}
\label{polc1=c4}
\prod_{\rho\in G} ((x+y)(y+z)(z+x)+\rho( c_1) (x^2 + y^2 + z^2+x y+   x z+  z y)+\rho(b_1)(x+y+z)+\rho(d)).
\end{equation}
The terms of degree 6 in the above expression are
\begin{eqnarray*}
&&((x+y)(y+z)(z+x))^2 (d+ \rho (d)+ \rho^2(d)) , \\
&&(x+y)(y+z)(z+x)  (x^2 + y^2 + z^2+x y+   x z+  z y)(x+y+z)\times\\
&&\qquad (c_1\rho(b_1) + b_1\rho(c_1) + c_1\rho^2(b_1) +\rho(c_1)\rho^2(b_1) + b_1\rho^2(c_1) +\rho(b_1)\rho^2(c_1)) ,\\
&&(x^2 + y^2 + z^2+x y+   x z+  z y)^3(c_1 \rho (c_1)\rho^2(c_1) ).
\end{eqnarray*}
Their sum must be equal to $\phi_9$.
We see \eb{by comparing the various monomials} that it implies that
$$d+ \rho (d)+ \rho^2(d)=c_1 \rho (c_1)\rho^2(c_1) $$
and
$$c_1\rho(b_1) + b_1\rho(c_1) + c_1\rho^2(b_1) +\rho(c_1)\rho^2(b_1) + b_1\rho^2(c_1) +\rho(b_1)\rho^2(c_1)=0.$$
The terms of degree 5 in the above expression are
\begin{eqnarray*}
&&(x+y)(y+z)(z+x)  (x^2 + y^2 + z^2+x y+   x z+  z y)\times\\
&&\qquad (c_1\rho(d_1) + d_1\rho(c_1) + c_1\rho^2(d_1) +\rho(c_1)\rho^2(d_1) + d_1\rho^2(c_1) +\rho(d_1)\rho^2(c_1)), \\
&&(x+y)(y+z)(z+x) (x+y+z)^2(  \rho(b_1)   \rho^2(b_1)+ \rho^2(b_1)    b_1+  b_1   \rho(b_1)) , \\
&& (x^2 + y^2 + z^2+x y+   x z+  z y)^2(x+y+z)\times\\
&&\qquad (c_1\rho(c_1)\rho^2(b_1)+\rho(c_1)\rho^2(c_1) b_1+\rho^2(c_1) c_1 \rho (b_1)).
\end{eqnarray*}
Their sum must be equal to $\phi_8$ which is zero. So one has the equations
\begin{eqnarray*}
c_1\rho(d_1) + d_1\rho(c_1) + c_1\rho^2(d_1) +\rho(c_1)\rho^2(d_1) + d_1\rho^2(c_1) +\rho(d_1)\rho^2(c_1)&=&0 , \\
 \rho(b_1)   \rho^2(b_1)+ \rho^2(b_1)    b_1+  b_1   \rho(b_1)&=&0,  \\
c_1\rho(c_1)\rho^2(b_1)+\rho(c_1)\rho^2(c_1) b_1+\rho^2(c_1) c_1 \rho (b_1)&=&0.  \\
\end{eqnarray*}

\subsubsection{Case where $c_1\ne0$}

We thus get the equations:\[
\left\{
\begin{array}{lcc}
 b_1+\rho(b_1)+\rho^2(b_1) & =  & 0 , \\
c_1\rho(b_1) + b_1\rho(c_1) + c_1\rho^2(b_1) +\rho(c_1)\rho^2(b_1) + b_1\rho^2(c_1) +\rho(b_1)\rho^2(c_1)  & =  &0 ,  \\
c_1\rho(c_1)\rho^2(b_1)+\rho(c_1)\rho^2(c_1) b_1+\rho^2(c_1) c_1 \rho (b_1) & =  &   0.
\end{array}
\right.
\]
We can consider this system as a homogeneous system of linear equations with $ b_1 $, $ \rho (b_1) $ and $ \rho ^ 2 (b_1) $ as unknown. Its determinant is
\begin{eqnarray*}
&& c_1^2 \rho (c_1) + c_1 \rho (c_1^2) + c_1^2 \rho^2(c_1) + \rho (c_1^2 )\rho^2(c_1)+ c_1 \rho^2(c_1^2) + \rho (c_1) \rho^2(c_1^2)\\
&&\qquad =(\rho(c_1) + \rho^2(c_1)) (c_1 + \rho(c_1)) (c_1 + \rho^2(c_1))\\
&&\qquad =c_1   \rho(c_1)  \rho^2(c_1)
\end{eqnarray*}
thanks to the equation
$c_1+\rho (c_1)+\rho^2(c_1)=0$.
If $ c_1 \ne0 $, the determinant is not zero, and the linear homogeneous system admits only a trivial solution:
$$b_1=\rho(b_1)=\rho^2(b_1)=0.$$

\subsubsection{Case where $c_1=0$ and $d\ne0$}

Looking at the terms of degree 7, 4 and 1, we find the equations
\[
\left\{
\begin{array}{lcc}
b + \rho (b_1) + \rho^2(b_1) &=& 0,\\
 d \rho (d) \rho^2(b_1) + d \rho (b_1) \rho^2 (d) + b \rho (d) \rho^2 (d) &=& 0,\\
  d \rho (b_1) + b \rho (d) + d \rho^2(b_1) +
    \rho (d) \rho^2(b_1) + b \rho^2 (d) + \rho (b_1) \rho^2 (d) &=& 0,
\end{array}
\right.
\]
    whose determinant is
\begin{eqnarray*}
&& d^2 \rho (d) + d \rho (d^2) + d^2 \rho^2(d) + \rho (d^2 )\rho^2(d)+ d \rho^2(d^2) + \rho (d) \rho^2(d^2)\\
&&\qquad =(\rho(d) + \rho^2(d)) (d + \rho(d)) (d + \rho^2(d))\\
&&\qquad =d   \rho(d)  \rho^2(d)
\end{eqnarray*}
thanks to the equation
$d+\rho (d)+\rho^2(d)=c_1\rho(c_1)\rho^2(c_1)=0$.
If $ d \ne0 $, the determinant is not zero, and the linear homogeneous system admits only a trivial solution:
$$b_1=\rho(b_1)=\rho^2(b_1)=0.$$

\subsubsection{Case where $c_1= d=0$}

The terms of degree 3 in the polynomial (\ref{polc1=c4}) are
$(x+y+z)^3$ with coefficient $b_1 \rho(b_1) \rho^2(b_1)$. These  terms must be equal to $a_6\phi_6$ with $\phi_6=(x+y)(y+z)(z+x)$.
This implies that $b_1 \rho(b_1) \rho^2(b_1)=0$ hence $b_1=0$.

 \subsubsection{Term of degree 4}

The terms of degree 4 of
$$\prod_{\rho\in G} ((x+y)(y+z)(z+x)+\rho( c_1) (x^2 + y^2 + z^2+x y+   x z+  z y) +\rho(d))$$
are
$(c_1 \rho (c_1) \rho^2(d)+\rho (c_1)\rho^2(c_1) d+\rho^2(c_1)  c_1 \rho (d))(x^2 + y^2 + z^2+x y+   x z+  z y) ^2$.
The polynomial $\phi_7$ is not of this form, thus one has $a_7=0$ and
$$c_1 \rho (c_1) \rho^2(d)+\rho (c_1)\rho^2(c_1) d+\rho^2(c_1)  c_1 \rho (d)=0.$$

 \subsubsection{Linear system}

Finally we get the equations
$$
\left\{
\begin{array}{lcc}
c_1 \rho (c_1) \rho^2(d)+\rho (c_1)\rho^2(c_1) d+\rho^2(c_1)  c_1 \rho (d) &=&0,\\
(\rho (c_1) +c_1) \rho^2(d)+( \rho^2(c_1)+\rho (c_1))  d+(c_1+\rho^2(c_1)) \rho (d)&=&0,\\
d+ \rho (d)+ \rho^2(d)+c_1 \rho (c_1)\rho^2(c_1)  &=   & 0 .
\end{array}
\right.
$$
They form  a linear system with $d$, $ \rho (d)$ et $\rho^2(d)$ as unknown.
Its determinant is
\begin{eqnarray*}
&& c_1^2 \rho (c_1) + c_1 \rho (c_1^2) + c_1^2 \rho^2(c_1) + \rho (c_1^2 )\rho^2(c_1)+ c_1 \rho^2(c_1^2) + \rho (c_1) \rho^2(c_1^2)\\
&&\qquad =(\rho(c_1) + \rho^2(c_1)) (c_1 + \rho(c_1)) (c_1 + \rho^2(c_1))\\
&&\qquad =c_1   \rho(c_1)  \rho^2(c_1)
\end{eqnarray*}
thanks to the equation
$c_1+\rho (c_1)+\rho^2(c_1)=0$.

If $c_1\ne0$, the determinant is nonzero, and the equation has only one solution
       \begin{eqnarray*}
 d&=&  { c_1^3\rho(c_1)\rho^2(c_1)          \over(c_1 -\rho(c_1)) (c_1 -\rho^2(c_1))}
 = { c_1^3\rho(c_1)\rho^2(c_1)          \over\rho^2(c_1) \rho(c_1)}
 =  { c_1^3}\\
\end{eqnarray*}
since  $c_1$  fulfills the equation  $c_1+\rho (c_1)+\rho^2(c_1)=0$ and
$ c_1  \ne  \rho (c_1)$ and    $c_1 \ne \rho^2(c_1)$.
The other solutions may be deduced by   Galois theory.
Then we have shown that the equation of
the surface $\cj X$ is
\begin{eqnarray*}
0=\phi&=&(A+c_1(x^2 + y^2 + z^2+x y+   x z+  z y)+c_1^3)\times\\
&&\qquad(A+ c_1^q(x^2 + y^2 + z^2+x y+   x z+  z y)+c_1^{3q}))\times\\
&&\qquad\qquad(A+c_1^{q^2}(x^2 + y^2 + z^2+x y+   x z+  z y)+c_1^{3q^2}) \\
&=&A^3+\beta A\mu^2+ (A^2+\mu^3)\gamma +A(\gamma^2+\beta ^3)+\beta ^2\gamma\mu+\gamma^3
\end{eqnarray*}
with
$\beta =c^{q+q^2}+c^{1+q^2}+c^{1+q}$,
$\mu=x^2 + y^2 + z^2+x y+   x z+  z y$
and
$\gamma=c^{1+q+q^2}$.

If $c_1=0$, the only condition on $d$ is:
$d+ \rho (d)+ \rho^2(d)=0$.
 The equation of the surface $\cj X$ is
$$0=\phi=(A+d)(A+d^q)(A+d^{q^2})=A^3+A(d^{1+q}+d^{1+q^2}+d^{q+q^2})+d^{1+q+q^2}.$$
\subsection{Proof of Theorem \ref{d12}}

\subsubsection{The functions $f$ associated to $X$, case $c_1\ne0$}

The computation of the coefficients $ a_i $ shows that the function $ f $ from which $\phi$ originates is
$$f(x)=(x^4+\beta x^2+\gamma x)^3 +\beta ^2x^8+\beta \gamma^2x^4. $$
This function (apart the terms $\beta ^2x^8+\beta \gamma^2x^4$ which can be eliminated by Proposition \ref{p2}) is the cube of the linear polynomial
$$x^4+x^2(c^{1+q}+c^{1+q^2}+c^{q+q^2})+x c^{1+q+q^2}.$$
This linear polynomial is equal to
$$x(x+c)(x+c^q)(x+c^{q^2}).$$
It has only one zero in $ \f q $, hence it is bijective.
Therefore the function $ f $ is CCZ-equivalent to the APN function $ x ^ 3 $.

\subsubsection{The functions $f$ associated to $X$, case $c_1=0$}
For $c=0$,
the function $ f $ from which $\phi$ originates is
$$f(x)=x^{12}+x^6(d^{1+q}+d^{1+q^2}+d^{q+q^2})+x^3d^{1+q+q^2}.$$
This function is composed of the linear polynomial
$$x^4+x^2(d^{1+q}+d^{1+q^2}+d^{q+q^2})+x d^{1+q+q^2}$$
  by the Gold-function $ x ^ 3 $. The linear polynomial is  bijective as before.
Therefore the function $ f $ is again CCZ-equivalent to the APN function $ x ^ 3 $.

It ends the proof of Theorem \ref{d12}.

\end{document}